\documentclass[12pt]{article}
\author{Constantin-Nicolae Beli\footnote{Supported by by the Romanian IDEI project PCE − 2012-4-364 of the Ministry of National Education CNCS-UEFISCIDI}} 
\title{Common slot lemmas in characteristic 2} 
\date{}
\usepackage{amsfonts,amsmath,amssymb}

   \def\m{\lim}
   
\def\p{\partial}     
    \def\te{\theta}
   
 \def\({\overline}

\def\){\underline} \def\<{\cdot} \def\go{\mathfrak}
\def\>{~~~~~~~} \def\#{{\bf
Definition}} \def\*{\section} \def\be{\begin{equation}}
\def\ee{\end{equation}}

\def\ti{\times}  \def\oo{{\cal O}} 
   
 \def\ff{\dot{F}} \def\ooo{{\oo^\ti}} 
 \def\mo{{\rm mod}~}  
  \def\fs{\ff^2}  
\def\p{\go p}    
\def\*{\sharp}  \def\0{} 
 \def\1{^{-1}}  
 \def\[{\prec} \def\]{\succ} 
\def\bmat{\left(\begin{array}} \def\emat{\end{array}\right)} 
\def\ap{\cong}   
\def\N{{\rm N}}  
\def\la{\langle} \def\ra{\rangle} 
 \def\m2{~(\mo 2)} \def\no{\noindent}
 \def\btm{\begin{thm}}
\def\etm{\end{tm}}
 \def\blem{\begin{lem}}
\def\elem{\end{lem}}

\newtheorem{theorem}{Theorem}[section]
\newtheorem{proposition}[theorem]{Proposition}
\newtheorem{lemma}[theorem]{Lemma}
\newtheorem{definition}{Definition}
\newtheorem{corollary}[theorem]{Corollary}

\newtheorem{bof}[theorem]{}
\newtheorem{teorema}{Theorem}

\def\qed{\mbox{$\Box$}\vspace{\baselineskip}}
\def\pf{$Proof.$} 
\def\bco{\begin{corollary}} \def\eco{\end{corollary}} 
\def\bdf{\begin{definition}} \def\edf{\end{definition}} 
\def\btm{\begin{theorem}} \def\etm{\end{theorem}} 
\def\bpr{\begin{proposition}} \def\epr{\end{proposition}}  
\def\blm{\begin{lemma}} \def\elm{\end{lemma}} 
\def\bff{\begin{bof}\rm} \def\eff{\end{bof}}
\def\btr{\begin{teorema}} \def\etr{\end{teorema}}

\def\de{\newcommand} \de\tm[1]{{\no\bf Theorem~#1}}

\de\lm[1]{{\no\bf Lemma~#1}}
\de\df[1]{{\no\bf Definition~#1}} \de\co[1]{{\no\bf Corollary~#1}}
\de\tp[1]{\te (#1 )} \de\ts[1]{\te (O^-(#1 ))} \de\ty[1]{\te
(O(#1 ))} \de\tx[1]{\te (#1 )} \de\up[1]{(1+\p^{#1} )\fs}
 \de\upn[2]{(1+\p^{#1})\fs\cap\N (#2 )} \de\xt[2]{\te (#1 /#2 )}
 \de\ups[1]{((1+\p^{#1})\fs )^*} \de\upo[1]{(1+\p^{#1} )\ooo^2}
\de\upon[2]{(1+\p^{#1})\ooo^2\cap\N (#2 )}
\de\lr[1]{\longrightarrow^{\!\!\!\!\!\!\!\! #1}}
\de\lf[1]{\longleftarrow^{\!\!\!\!\!\!\!\! #1}}
\de\si[1]{\sim^{\!\!\!\!\! #1}} \de\apr[1]{\approx^{\!\!\!\!\! #1}}
\de\leg[2]{\left(\frac {#1}{#2}\right)}
\DeclareMathOperator\Br{Br}

\def\Brd{{}_2\Br}

\DeclareMathOperator\car{char}

\begin{document}

\maketitle

\begin{abstract}

The Hilbert symbol $(\cdot,\cdot )$ from characteristic $\neq 2$ has
two analogues in characteristic $2$, $[\cdot,\cdot )$ and
$((\cdot,\cdot ))$. 

The well known common slot lemma from characteristic $\neq 2$ (see
e.g. [T, Theorem 4.4 and the following Remark 3)]), which involves the
symbol $(\cdot,\cdot )$ has three analogues in characteristic $2$. Two
of these analogues involve the symbol $[\cdot,\cdot )$ and one
involves $((\cdot,\cdot ))$. Of these three analogues so far only
one has been considered. (See [D, \S 14, Theorem 7], where it is
stated as an exercise.) In this paper we state and prove the
remaining two analogues. 



\end{abstract}



\section{The symbols $(\cdot,\cdot )$, $[\cdot,\cdot )$ and
$((\cdot,\cdot ))$. Definition and basic properties}

If $A$ is a central division algebra over a field $F$ then we denote
by $[A]$ its class in the Brauer group $(\Br (F),+)$. We have
$[A]+[B]=[A\otimes B]$. We also denote by $\Brd (F)$ the 2-torsion of
$\Br (F)$. 

We denote by $\la a_1,\ldots,a_n\ra$ the diagonal quadratic form
$a_1X_1^2+\cdots +a_nX_n^2$, and by $[a,b]$ the binary quadratic form
$aX^2+XY+bY^2$.

If $\car F\neq 2$ then we have the Hilbert symbol
$$(\cdot,\cdot ):\ff/\fs\ti\ff/\fs\to\Brd (F)$$. 
If $a,b\in\ff$ then $(a,b):=[Q_{(a,b)}]$, where $Q_{(a,b)}$ is the
quaternion algebra generated by $1,i,j,ij$, with the relations
$i^2=a$, $j^2=b$, $ij+ji=0$. The norm map $N:Q_{(a,b)}\to F$ is given
by $Q_{(a,b)}(X+Yi+Zj+Tij)=X^2-aY^2-bZ^2+abT^2$ so
$(Q_{(a,b)},N)\ap\la 1,-a,-b,ab\ra$ relative to the basis $1,i,j,ij$.


If $\car F=2$ then, besides $(\ff/\fs,\cdot )$, we also have the
groups $(F/F^2,+)$ and $(F/\wp (F),+)$, where $\wp :F\to F$ is the
Artin-Schreier map, $\wp (x)=x^2+x$. Then we define two symbols:
$$[\cdot,\cdot ):F/\wp (F)\ti\ff/\fs\to\Brd (F),$$
is given by $[a,b)=[Q_{[a,b)}]$, where $Q_{[a,b)}=\la 1,i,j,ij\ra$, with
the relations $i^2+i=a$, $j^2=b$ and $ij+ji=j$. The norm map
$N:Q_{[a,b)}\to F$ is given by
$Q_{[a,b)}(X+Yi+Zj+Tij)=X^2+XY+aY^2+b(Z^2+ZT+aT^2)$ so 
$(Q_{(a,b)},N)\ap [1,a]\perp b[1,a]$ relative to the basis $1,i,j,ij$.


$$((\cdot,\cdot )):F/F^2\ti F/F^2\to\Brd (F),$$
is given by $((a,b))=[Q{((a,b))}]$, where $Q_{((a,b))}=\la
1,i,j,ij\ra$ with the relations $i^2=a$, $j^2=b$ and $ij+ji=1$. The
norm map $N:Q_{((a,b))}\to F$ is given by
$Q_{(a,b)}(X+Yi+Zj+Tij)=X^2+XT+abT^2+aY^2+YZ+bZ^2$ so 
$(Q_{((a,b))},N)\ap\la [1,ab]\perp [a,b]$ relative to the basis
$1,ij,i,j$.


The symbols $(\cdot,\cdot )$, $[\cdot,\cdot )$ and $((\cdot,\cdot ))$
are bilinear, $(\cdot,\cdot )$ is symmetric and $((\cdot,\cdot ))$ 
is anti-symmetric. In particular, because of the $2$-torsion,
$((\cdot,\cdot ))$ is also symmetric. (We have $((a,b))=-((b,a))=((b,a))$.) 

The symbols $[\cdot,\cdot )$ and $((\cdot,\cdot ))$ are related by the
relations $((a,b))=[ab,b)$ if $b\neq 0$ and $[a,0)=0$. 
\bigskip

Note that, unlike $(\cdot,\cdot )$ and $[\cdot,\cdot )$, the symbol
$((\cdot,\cdot ))$ is not very well known and used. Most authors
simply write $[ab,b)$ instead of $((a,b))$. Also the notation
$((\cdot,\cdot ))$ is not universally accepted. For properties of
$((\cdot,\cdot ))$ see e.g. [EKM, \S 98.E], where $((a,b))$ is
denoted by $\left[\frac{a,b}F\right]$. 

\section{Common slot lemma}

In characteristic $\neq 2$ we have the following result known as
``chain lemma'' or ``common slot Lemma''.

\btm  Let $Q$ be a quaternion algebra over a field $F$ with $\car
F\neq 2$. If $[Q]=(a_1,b_1)=(a_2,b_2)$ then there is $b\in\ff$ such
that $[Q]=(a_1,b)=(a_2,b)$. 
\etm

The chain lemma has three analogues in characteristic $2$.

\btm Let $Q$ be a quaternion algebra over a field $F$ with $\car
F=2$.

(i) If $[Q]=[a_1,b_1)=[a_2,b_2)$ then there is $b\in\ff$ such
that $[Q]=[a_1,b)=[a_2,b)$. 

(ii) If $[Q]=[a_1,b_1)=[a_2,b_2)$ then there is $a\in F$ such that
$[Q]=[a,b_1)=[a,b_2)$. 

(iii) If $[Q]=((a_1,b_1))=((a_2,b_2))$ then there is $a\in F$ such that
$[Q]=((a,b_1))=((a,b_2))$. 
\etm

For the proof of these statements we use notations and results from
[SV, Chapter 1]. We denote by $N_Q:Q\to F$ the norm map of $Q$ and by
$b_Q$ the corresponding polar map, $b_Q(x,y)=N_Q(x+y)-N_Q(x)-N_Q(y)$.

In both Theorems 2.1 and 2.2 we may assume that $Q$ is non-split,
i.e. $(Q,N_Q)$ is unisotropic since otherwise $Q\ap M_2(F)$ so
$[Q]=0$. Then in Theorem 2.1 amd 2.2(i) we may take $b=1$ (we have
$(a_i,1)=0$ and $[a_i,1)=0$, respectively) and in Theorem 2.2(ii)
and (iii) we may take $a=0$ (we have $[0,b_i)=0$ and $((0,b_i))=0$,
respectively). 
\medskip

{\em Proof of Theorem 2.1 and Theorem 2.2(i).} For $i=1,2$ we
have $Q=\langle e,x_i,y_i,x_iy_i\rangle$, where in the case of (i)
$x_i^2=a_ie$, $y_i^2=b_ie$ and $x_iy_i+y_ix_i=0$ and in the case of
(ii) $x_i^2+x_i=a_ie$, $y_i^2=b_ie$ and $x_iy_i+y_ix_i=y_i$. If
$D_i=\langle e,x_i\rangle$ then $D_i$ is a degree two extension of
$k$, $D_i\cong k(\sqrt{a_i})$ in the case of (i) and $D_i\cong
k(\wp^{-1}(a_i))$ in the case of (ii). 

The othogonal complement of $\langle e,x_1,x_2\rangle$ is $\neq 0$ so
let $y\in\langle e,x_1,x_2\rangle^\perp$, $y\neq 0$. We have $y\in
D_i^\perp$ for $i=1,2$ and since $N_Q$ is anisotropic we have
$-b:=N_C(y)\neq 0$. Then $Q$ is obtained from $D_i$ by the doubling
process from [SV, Prposition 1.5.1], $Q=D_i\oplus D_iy$. Since
$D_i\cong k(\sqrt{a_i})$ or $k(\wp^{-1}(a_i))$ and $N_Q(y)=-b$ we have
$[Q]=(a_i,b)$ or $[a_i,b)$, respectively, for $i=1,2$. \qed

For Theorem 2.2(ii) and (iii) we need the following result.

\blm Let $C$ be a composition algebra over a field $F$ of
characteristic $2$ and let $b_1,b_2\in k$. If there are $y_1,y_2\in
e^\perp\setminus Fe$ with $N_C(y_i)=b_i$ then such $y_1,y_2$ can be
chosen with the additional property that $y_1\not\perp y_2$.
\elm
\pf  Assume that $y_1\perp y_2$. First we prove that there is some
$u\in C$ with $u\perp e$ but $u\not\perp y_1,y_2$. If $e,y_1,y_2$ are
not linearly independent then $y_2=\alpha e+\beta y_1$ for some
$\alpha,\beta\in F$, $\beta\neq 0$. (We have $y_1,y_2\notin Fe$.) By
the regularity of $b_C$ there is some $u\in C$ with $b_C(u,e)=0$ and
$b_C(u,y_1)=1$. Then we also have $b(u,y_2)=\beta\neq 0$ so $u\perp e$
and $u\not\perp y_1,y_2$. If $e,y_1,y_2$ are linearly independent
then, again by the regularity of $b_C$, there is $u\in C$ with
$b_C(u,e)=0$ and $b_C(u,y_1)=b_C(u,y_2)=1$ and again $u\perp e$ and
$u\not\perp y_1,y_2$. If $N_C(u)=0$ then let $u'=u+e$. Since $e\perp
e,y_1,y_2$ and $u\perp e$, $u\not\perp y_1,y_2$ we have $u'\perp e$,
$u'\not\perp y_1,y_2$. However $N_C(u')=N_C(u)+N_C(e)=1\neq 0$. So we
may assume that $N_C(u)\neq 0$ and we may consider the orthogonal
transvection $\sigma\in O(N_C)$ given by $x\mapsto
x+b_c(u,x)N_C(u)^{-1}u$. Let $y'_2=\sigma
(y_2)=y_2+b_c(u,y_2)N_C(u)^{-1}u$. We have $N_C(y'_2)=N_C(y_2)=b_2$,
$b_C(y'_2,e)=b_C(y_2,e)+b_c(u,y_2)N_C(u)^{-1}b_C(u,e)=0$ and
$b_C(y'_2,y_1)=b_C(y_2,y_1)+b_c(u,y_2)N_C(u)^{-1}b_C(u,y_1)=
b_c(u,y_2)N_C(u)^{-1}b_C(u,y_1)\neq 0$. So $y'_2\perp e$ but
$y'_2\not\perp y_1$. Thus $y_1,y'_2$ satisfy all the required
conditions. \qed
\medskip

{\em Proof of Theorem 2.2(ii) and (iii).} As seen from $\S$1, for
$i=1,2$, in the case of (i) we have $(Q,N_Q)\cong [1,a_i]\perp
b_i[1,a_i]$ relative to some basis $e,x_i,y_i,x_iy_i$; and in the case
of (ii) we have $N_Q\cong [1,a_ib_i]\perp [a_i,b_i]$ relative to some basis
$e,x_iy_i,x_i,y_i$. In both cases $N_Q(y_i)=b_i$ and $e\perp y_i$ for
$i=1,2$. We also have $y_i\notin Fe$. Then we can use Lemma 2 and we
can choose $y_1,y_2\in e^\perp\setminus Fe$ such that $N_Q(y_i)=b_i$
with the additional property that $y_1\not\perp y_2$. We cannot have
$y_2\in\langle e,y_1\rangle$ for this would imply that $y_2\perp
y_1$. Thus $e,y_1,y_2$ are linear independent. We now treat separately
the statements (ii) and (iii) of Theoren 2.2.

(ii) By the regularity of $b_Q$ there is some $x\in Q$ with
$b_Q(x,e)=1$ and $b_C(x,y_1)=b_C(x,y_2)=0$. Let $a=N_Q(x)$. Since
$b_Q(x,e)=1$, by [SV, Proposition 1.2.3] we have $x^2+x+ae=0$. Thus
$D=\langle e,x\rangle$ is a field extension of $F$, $D\cong
F(\wp^{-1}(a))$. Since $y_i\perp e,x$, so $y_i\in D^\perp$, and
$N_Q(y_i)=b_i\neq 0$ we have by the doubling process of [SV,
Proposition 1.5.1] $Q=D\oplus Dy_i$ and $[Q]=[a,b_i)$. 

(iii) By [SV, Proposition 1.2.3] we have $y_i^2+b_ie=0$,
i.e. $y_i^2=b_ie$. Also $\overline{y_i}=y_i$. By the regularity of
$b_Q$ there is some $x\in Q$ with $b_Q(x,e)=0$ and
$b_C(x,y_1)=b_C(x,y_2)=1$. By [SV, (1.3)] we have
$b_Q(xy_i,y_i)=b(x,e)N_Q(y_i)=0$ and by [SV, Lemma 1.3.2] we have
$b_Q(xy_i,e)=b_Q(x,\overline{y_i})=b_Q(x,y_i)=1$. So $e,xy_i,y_i$
are in the same situation as $e,x,y_i$ were in the proof of (ii). Since
$N_Q(xy_i)=N_Q(x)N_Q(y_i)=ab_i$ and $N_Q(y_i)=b_i$ we have
$[Q]=[ab_i,b_i)=((a_i,b_i))$. \qed
\bigskip

{\bf Remark 1} The proof of the common slot lemma in characteristic
$\neq 2$ (our Theorem 2.1) which we provided here is essentially the same
as the one from [T, Theorem 4.4 and the following Remark 3)]. The fact
that the proof of Theorem 2.2(i) goes on the same lines shows that
Theorem 2.2(i) is the ``correct'' analogue of the common slot theorem
in characteristic $2$. The result stated in [D, \S 14, Theorem 7] is
our Theorem 2.2(ii). The author leaves this result as an exercise,
pointing to Tate's paper for the characteristic $\neq 2$ case. Perhaps
by mistake he stated the ``wrong'' analogue, which turns out to be
true but by a proof that is not similar to the one from
characteristic $\neq 2$. 

{\bf Remark 2} One may prove the common slot lemma in characteristic
$\neq 2$ by using only the quadratic space structure on $Q$, and
ignoring the algebraic structure. That is, we use the information that
$(Q,N_Q)\ap\la 1, -a_1,-b_1,a_1b_1\ra\ap \la 1, -a_2,-b_2,a_2b_2\ra$
to prove that there is some $b\in\ff$ such that $(Q,N_Q)\ap\la 1,
-a_1,-b,a_1b\ra\ap \la 1, -a_2,-b,a_2b\ra$. This works very nice in
characteristic $\neq 2$. We can do the same for Theorem 2.2 but the
proof would be longer and messy. 
\bigskip

{\bf References}

[D] P. K. Draxl, ``Skew Fields'', London Mathematical Society Lecture
Note Series. 81. 

[EKM] Richard Elman, Nikita Karpenko and Alexander Merkurjev, ``The
Algebraic and Geometric Theory of Quadratic Forms'', AMS Colloquim
Publications, vol. 56. 

[SV] Tonny A. Springer and Ferdinand D. Veldkamp, ``Octonions, Jordan Algebras
and Exceptional Groups'', Springer Monographs in Mathematics, 2000. 

[T] John Tate, ``Relations between $K_2$ and Galois Cohomology'',
Inventiones math. 36, 257-274 (1976).
\bigskip

Institute of Mathematics Simion Stoilow of the Romanian Academy, Calea Grivitei 21, RO-010702 Bucharest, Romania.

E-mail address: Constantin.Beli@imar.ro

\end{document}